\documentclass[11pt]{article}
\usepackage[english]{babel}
\usepackage[latin1]{inputenc}
\usepackage{exscale}
\usepackage[centertags]{amsmath}
\usepackage{amsfonts,amstext,amssymb,amsthm}
\usepackage{comment}
\usepackage{newlfont}
\usepackage{enumerate}
\usepackage{graphicx}
\usepackage{ContinuidadeAbraham7}
\usepackage{color}
\usepackage[usenames,dvipsnames]{xcolor}
\definecolor{BlueFonse}{rgb}{0,0,1}
\definecolor{BlueFonse1}{cmyk}{1,0,0,0.7}
\usepackage[colorlinks=true,linkcolor=blue,citecolor=ForestGreen]{hyperref}
\usepackage{hyperref}
\usepackage{url}
\title{Continuity and differentiability properties of the isoperimetric profile in complete noncompact Riemannian manifolds with bounded geometry}
\author{Abraham Mu\~noz Flores\footnote{Partially supported by Capes}, Stefano Nardulli}
\begin{document}
      \maketitle
\pretolerance=2000
\tolerance=3000
\noindent {\sc abstract}. For a complete noncompact connected Riemannian manifold with bounded geometry $M^n$, we prove that the isoperimetric profile function $I_{M^n}$ is continuous. Here for bounded geometry we mean that $M$ have $Ricci$ curvature bounded below and volume of balls of radius $1$, uniformly bounded below with respect to its centers. Then under an extra hypothesis on the geometry of $M$, we apply this result to prove some differentiability property of $I_M$ and a differential inequality satisfied by $I_M$, extending in this way well known results for compact manifolds, to this class of noncompact complete Riemannian manifolds with bounded geometry.
\bigskip\bigskip

\bigskip

\centerline{\bf AMS subject classification: }
49Q20, 58E99, 53A10, 49Q05.
      \tableofcontents       
\section{Introduction}\label{1} 
In the remaining part of this paper we always assume that all the Riemannian manifolds $(M, g)$ considered are smooth with smooth Riemannian metric $g$. We denote by $V$ the canonical Riemannian measure induced on $M$ by $g$, and by $A$ the $(n-1)$-Hausdorff measure associated to the canonical Riemannian length space metric $d$ of $M$. When it is already clear from the context, explicit mention of the metric $g$ will be suppressed in what follows. We give here the basic definitions of $BV$-functions and finite perimeter sets on a manifold. 
\begin{Def} Let $M$ be a Riemannian manifold of dimension $n$, $U\subseteq M$ an open subset, $\mathfrak{X}_c(U)$ the set of smooth vector fields with compact support on $U$. Given a function $u\in L^{1}(M)$, define the variation of $u$ by
      \begin{equation}
                 |Du|(M):=sup\left\{\int_{M}u div_g(X)dV_g: X\in\mathfrak{X}_c(M), ||X||_{\infty}\leq 1\right\},
      \end{equation}  
      where $||X||_{\infty}:=\sup\left\{|X_p|_{{g}_p}: p\in M\right\}$ and $|X_p|_{{g}_p}$ is the norm of the vector $X_p$ in the metric $g_p$ on $T_pM$.
      We say that a function $u\in L^{1}(M)$, has \textbf{bounded variation}, if $|Du|(M)<\infty$ and we define the set of all functions of bounded variations on $M$ by $BV(M):=\{u\in L^1(M):\:|Du|(M)<+\infty\}$.
      \end{Def}
\begin{Def} Let $M$ be a Riemannian manifold of dimension $n$. Given $E\subset M$ measurable with respect to the Riemannian measure, $U\subseteq M$ an open subset, the \textbf{perimeter of $E$ in $U$}, $ \mathcal{P}(E, U)\in [0,+\infty]$, is defined as
      \begin{equation}
                 \mathcal{P}(E, U):=sup\left\{\int_{U}\chi_E div_g(X)dV_g: X\in\mathfrak{X}_c(U), ||X||_{\infty}\leq 1\right\},
      \end{equation}  
where $||X||_{\infty}:=\sup\left\{|X_p|_{{g}_p}: p\in M\right\}$ and $|X_p|_{{g}_p}$ is the norm of the vector $X_p$ in the metric $g_p$ on $T_pM$. If $\mathcal{P}(E, U)<+\infty$ for every open set $U$, we call $E$ a \textbf{locally finite perimeter set}. Let us set $\mathcal{P}(E):=\mathcal{P}(E, M)$. Finally, if $\mathcal{P}(E)<+\infty$ we say that \textbf{$E$ is a set of finite perimeter}.    
\end{Def}
When dealing with finite perimeter sets or locally finite perimeter sets we will denote the reduced boundary $\partial^*\Omega$, by $\partial\Omega$ when no confusion may arise. For this reason we will denote $\mathcal{P}(\Omega)=A(\partial^*\Omega)=A(\partial\Omega)$ and for every finite perimeter set $\Omega'$ we always choose a representative $\Omega$ (i.e., that differs from $\Omega'$ by a set of Riemannian measure $0$), such that $\partial_{top}\Omega=\overline{\partial^*\Omega}$, where $\partial_{top}\Omega$ is the topological boundary of $\Omega$. At this point we give the definition of the isoperimetric profile function which is the main object of study in this paper.   
\subsection{The isoperimetric profile}
\begin{Def}\label{Def:IsPStrong}
Typically in the literature, the \textbf{isoperimetric profile function of $M$} (or briefly, the isoperimetric profile) $I_M:[0,V(M)[\rightarrow [0,+\infty [$, is defined by $$I_M(v):= \inf\{A(\partial \Omega): \Omega\in \tau_M, V(\Omega )=v \},$$ where $\tau_M$ denotes the set of relatively compact open subsets of $M$ with smooth boundary.
\end{Def}
However there is a more general context in which to consider this notion that will be better suited to our purposes.  Namely, we can give a weak formulation of the preceding variational problem replacing the set $\tau_M$ with the family of subsets of finite perimeter of $M$. 
\begin{Def}\label{Def:IsPWeak}
Let $M$ be a Riemannian manifold of dimension $n$ (possibly with infinite volume). We denote by $\tilde{\tau}_M$ the set of  finite perimeter subsets of $M$. The function $\tilde{I}_M:[0,V(M)[\rightarrow [0,+\infty [$  defined by 
     $$\tilde{I}_M(v):= \inf\{\mathcal{P}(\Omega): \Omega\in \tilde{\tau}_M, V(\Omega )=v\},$$ 
is called the \textbf{weak isoperimetric profile function} (or shortly the \textbf{isoperimetric profile}) of the manifold $M$. If there exists a finite perimeter set $\Omega\in\tilde{\tau}_M$ satisfying $V(\Omega)=v$, $\tilde{I}_M(V(\Omega))=A(\partial^*\Omega)= \mathcal{P}(\Omega)$ such an $\Omega$ will be called an \textbf{isoperimetric region}, and we say that $\tilde{I}_M(v)$ is \textbf{achieved}. 
\end{Def} 
There are many others possible definitions of isoperimetric profile corresponding to the minimization over various different admissible sets, as stated in the following definition.  
\begin{Def} For every $v\in[0, +\infty[$, let us define
\begin{eqnarray*}
I^{*}_{M}(v):=inf\{A(\partial_{top}\Omega): \Omega\subset M, \partial_{top}\Omega\: \text{is}\: C^{\infty}, V(\Omega)=v\},\\
\tilde{I}_{M}^{*}(v):=inf\{\mathcal{P}_{M}(\Omega):\Omega\subset M, \Omega\in\tilde{\tau}_{M} , V(\Omega)=v, diam(\Omega)<+\infty\}, 
\end{eqnarray*}
where $diam(\Omega):=\sup\{d(x,y):x,y\in\Omega\}$ denotes the \textbf{diameter} of $\Omega$.
\end{Def}
\begin{Rem}\label{Rem:Trivialinequality}
Trivially one have $I_M\geq I^{*}_{M}\geq\tilde{I}_M$ and $I_M\geq\tilde{I}^*_{M}\geq\tilde{I}_M$. 
\end{Rem}
However as we will see in Theorem $\ref{Thm:Equivalence}$, all of these definitions are actually equivalents, in the sense that the infimum remains unchanged, i.e., $I_M=\tilde{I}_M$. The proof of this fact involves actually very natural ideas. In spite of this it is technical and we have found no written traces in the literature, unless Lemma $2$ of \cite{Modica} that deal with the case of a compact domain of $\R^n$ as an ambient space. Hence we provided ourselves a proof. This equivalence allows us to consider elements of $\tau_M$ or $\tilde{\tau}_M$ according to what is more convenient to us. This observation is used in a crucial way when we prove Theorem $\ref{Main}$, see for example the proof of inequality \eqref{Eq:Main1}.    
The next fact to be observed is that it is worth to have a proof of the continuity of the isoperimetric profile, because in general the isoperimetric profile function of a complete Riemannian manifold is not continuous. In case of manifolds with density, in Proposition $2$ of \cite{MorganBlog}  is exhibited an example of  a manifold with density having discontinuous isoperimetric profile. To exhibit a complete Riemannian manifold with a discontinuous isoperimetric profile is a more  subtle and difficult task that was performed by the second author and Pierre Pansu in \cite{NardulliPansuDiscontinuous}, for manifolds of dimension $n\geq 3$, but whose methods with a slight modification of the arguments could be used also to settle the case $n=2$.  In spite of these quite sophisticated counterexamples the class of manifolds admitting a continuous isoperimetric profile is vast, for an account of the existing literature on the continuity results obtained for $I_M$, one could consult the introduction of \cite{RitoreContinuity} and the references therein. If $M$ is compact, classical compactness arguments of geometric measure theory combined with the direct method of the calculus of variations provide a short proof of the continuity of $I_{M}$ in any dimension $n$, \cite{MorganBlog} Proposition $1$. Finally, if $M$ is complete, non-compact, and $V(M)<+\infty$, an easy consequence of Theorem $2.1$ in \cite{RRosales} yields the possibility of extending the same compactness argument valid in the compact case and to prove the continuity of the isoperimetric profile, see for instance Corollary 2.4 of \cite{NardulliRusso}. A careful analysis of Theorem $1$ of \cite{Nar12} about the existence of generalized isoperimetric regions, leads to the continuity of the isoperimetric profile $I_M$ in manifolds with bounded geometry satisfying some other assumptions on the geometry of the manifold at infinity, of the kind considered by the second author and A. Mondino in \cite{MonNar}, i.e., for every sequence of points diverging to infinity, there exists a pointed smooth manifold $(M_{\infty}, g_{\infty}, p_{\infty})$ such that $(M,g,p_j)\rightarrow (M_{\infty}, g_{\infty}, p_{\infty})$ in $C^0$-topology. This proof is independent from that of Theorem \ref{Main}. This is not the case for general complete infinite-volume manifolds $M$. Recently Manuel Ritor\'e (see for instance \cite{RitoreContinuity}) showed that a complete Riemannian manifold possessing a strictly convex Lipschitz continuous exhaustion function has continuous  and nondecreasing isoperimetric profile $\tilde{I}_M$. Particular cases of these manifolds are Cartan-Hadamard manifolds and complete noncompact manifolds with strictly positive sectional curvatures. In \cite{RitoreContinuity} as in our Theorem \ref{Main} the major difficulty consists in finding a suitable way of subtracting a volume to an almost minimizing region. 

The aim of this paper is to prove Theorem \ref{Main} in which we give a very short and quite elementary proof of the continuity of $I_M$ when $M$ is a complete noncompact Riemannian manifold of bounded geometry. The reason which allow us to achieve this goal, is that in bounded geometry it is always possible to add or subtract to a finite perimeter set a small ball that captures a fixed fraction of volume (depending only from the bounds of the geometry) centered at points close to it. Following this philosophy it is quite easy to show that to have an isoperimetric region of volume $v$ ensures the upper semicontinuity of $I_M$ at $v$. This is exactly the content of Theorem \ref{Thm:Uppersemicontinuity}, in which we are also more lucky and we can subtract a ball of the right volume entirely contained in the isoperimetric region. The problems appears when we try to prove lower semicontinuity. To prove lower semicontinuity we need some kind of compactness that is expressed here by a bounded geometry condition.  Geometrically speaking our assumptions of bounded geometry ensures that the manifold at infinity is not too thin and enough thick to permit to place a small geodesic ball $B$ close to an arbitrary domain $D$ in such a way $V(B\cap D)$ recovers a controlled fraction of $V(D)$ and this fraction depends only on $V(D)$ and the bounds on the geometry $n, v_0, k$, see Definition 
\ref{Def:BoundedGeometry} below for the exact meaning of $n$, $v_0$, $k$.  The proof that we present here uses only metric properties of the manifolds with bounded geometry and for this reason it is still valid when suitably reformulated in the context of metric measure spaces. One can find similar ideas alredy in the metric proof of continuity of the isoperimetric profile contained in \cite{Gallot}. For the full generality of the results we need that the spaces have to be doubling, satisfying a $1$-Poincar\'e inequality and a curvature dimension condition. This class of metric spaces includes for example manifolds with density as well  as subRiemannian manifolds.  
We observe that another proof of Corollary \ref{BPGen} is possible following the same lines of \cite{BP}, the arguments used there permits also to obtain another proof of the continuity of the isoperimetric profile under our assumptions of bounded geometry but with the extra assumption of the existence of isoperimetric regions of every volume, which is less general of our own proof of Theorem $\ref{Main}$, because in Theorem $\ref{Main}$ we do not need to assume any kind of existence of isoperimetric regions. In spite of this the Heintze-Karcher type arguments used in \cite{BP} have an advantage because they permits to give a uniform bound on the length of the mean curvature vector of the generalized isoperimetric regions (i.e., left and right derivatives of $I_M$) with volumes inside an interval $[a, b]\subset ]0, V(M)[$, depending only on $a$ and $b$. Finally, we mention that just with Ricci bounded below and existence of isoperimetric regions the arguments of \cite{BP} fails and we cannot prove the continuity of the isoperimetric profile, for this we need a noncollapsing condition on the volume of geodesic balls as in our definition of bounded geometry. We give a detailed account of these arguments in Theorem \ref{Thm:BavardPansu}.
\subsection{Plan of the article}
\begin{enumerate}
           \item  Section $\ref{1}$ constitutes the introduction of the paper. We state the main results of the paper.
           \item In Section $\ref{Sec:Weak}$ we prove that $\tilde{I}_M=I_M$.
           \item In section $\ref{Sec:Continuity}$ we prove the continuity of isoperimetric profile in bounded geometry, i.e., Theorem $\ref{Main}$, without assuming existence of isoperimetric regions. 
           \item In  section $\ref{Sec:Differentiability}$, we prove Corollary $\ref{BPGen}$ and $\ref{MJGen}$. 
           \item In section $\ref{Apendix:BP}$ we explain the link with the preexisting work \cite{BP}.
\end{enumerate}
\subsection{Acknowledgements}  
The second author is indebted to Pierre Pansu, Frank Morgan, Andrea Mondino, and Luigi Ambrosio for useful discussions on the topics of this article. 
The first author wish to thank the CAPES for financial support.   
\subsection{Main Results}\label{MainRes} 
\begin{Res}\label{Thm:Equivalence} If $M^n$ is an arbitrary complete Riemannian manifold, then $I_M(v)=\tilde{I}_{M}^{*}(v)=\tilde{I}_M(v)=I^{*}_{M}(v)$. 
\end{Res} 
\begin{Def}\label{Def:BoundedGeometry}
A complete Riemannian manifold $(M, g)$, is said to have \textbf{bounded geometry} if there exists a constant $k\in\mathbb{R}$, such that $Ric_M\geq k(n-1)$ (i.e., $Ric_M\geq (n-1)kg$ in the sense of quadratic forms) and $V(B_{(M,g)}(p,1))\geq v_0$ for some positive constant $v_0$, where $B_{(M,g)}(p,r)$ is the geodesic ball (or equivalently the metric ball) of $M$ centered at $p$ and of radius $r> 0$.
\end{Def}
\begin{Res}[Continuity of the isoperimetric profile]\label{Main}
Let $M^n$ be a complete smooth Riemannian manifold with $Ric_M\geq (n-1)k$, $k\in\mathbb{R}$ and $V(B(p,1))\geq v_0>0$. Then $I_M$ is continuous on $[0,V(M)[$.
\end{Res}
\begin{Def} For any $m\in\mathbb{N}$, $\alpha\in [0, 1]$, a sequence of pointed smooth complete Riemannian manifolds is said to \textbf{converge in
the pointed $C^{m,\alpha}$, respectively $C^{m}$ topology to a smooth manifold $M$} (denoted $(M_i, p_i, g_i)\rightarrow (M,p,g)$), if for every $R > 0$ we can find a domain $\Omega_R$ with $B(p,R)\subseteq\Omega_R\subseteq M$, a natural number $\nu_R\in\mathbb{N}$, and $C^{m+1}$ embeddings $F_{i,R}:\Omega_R\rightarrow M_i$, for large $i\geq\nu_R$ such that $B(p_i,R)\subseteq F_{i,R} (\Omega_R)$ and $F_{i,R}^*(g_i)\rightarrow g$ on $\Omega_R$ in the $C^{m,\alpha}$, respectively $C^m$ topology. 
\end{Def}
\begin{Def}\label{Def:BoundedGeometryInfinity}
We say that a smooth Riemannian manifold $(M^n, g)$ has $C^{m,\alpha}$-\textbf{locally asymptotic bounded geometry} if it is of bounded geometry and if for every diverging sequence of points $(p_j)$, there exists a subsequence $(p_{{j}_{l}})$ and a pointed smooth manifold $(M_{\infty}, g_{\infty}, p_{\infty})$ with $g_{\infty}$ of class $C^{m,\alpha}$ such that the sequence of pointed manifolds $(M, p_{{j}_{l}}, g)\rightarrow (M_{\infty}, g_{\infty}, p_{\infty})$, in  $C^{m,\alpha}$-topology.  
\end{Def}
\begin{CorRes}[Bavard-Pansu-Morgan-Johnson in bounded geometry]\label{BPGen}  Let $M$ have $C^0$-locally asymptotic bounded geometry in the sense of Definition \ref{Def:BoundedGeometryInfinity}. Suppose that all the limit manifolds have a metric at least of class $C^2$. Then $I_M$ is absolutely continuous and twice differentiable almost everywhere. The left and right derivatives $I_M^-\geq I_M^+$ exist everywhere and their singular parts are non-increasing. If $k>0$ then $I_M$ is strictly concave on $]0, V(M)[$. If $k=0$, then $I_M$ is just concave on $]0, V(M)[$. If $k<0$, then $I_M(v)+C(a,b)v^2$ is concave, ($I_M$ could not be concave). Moreover, we have for every $k\in\R$ and almost everywhere
\begin{equation}\label{GenBP}
          I_MI_M^{''}\leq-\frac{{I'}_M^2}{n-1}-(n-1)k,
\end{equation}
with equality in the case of the simply connected space form of constant sectional curvature $k$. In this case, a generalized isoperimetric region is totally umbilic. 
\end{CorRes}
\begin{CorRes}[Morgan-Johnson isoperimetric inequality in bounded geometry]\label{MJGen} Let $M$ have $C^{2,\alpha}$-bounded geometry, sectional curvature $K$ and Gauss-Bonnet-Chern integrand $G$. Suppose that
\begin{itemize}
          \item $K<K_0$, or
          \item $K\leq K_0$, and $G\leq G_0$,
\end{itemize}
where $G_0$ is the Gauss-Bonnet-Chern integrand of the model space form of constant curvature $K_0$. Then for small prescribed volume, the area of a region $R$ of volume $v$ is at least as great as $A(\partial B_v)$, where $B_v$ is a geodesic ball of volume $v$ in the model space, with equality only if $R$ is isometric to $B_v$.   
\end{CorRes}
The proofs of Corollaries \ref{BPGen} and \ref{MJGen} run along the same lines as the corresponding proofs of Theorems $3.3$ and $4.4$ of \cite{MJ}. These last theorem where proven for compact manifolds the  needed changes to make them works in our more general context are fully provided in Section $\ref{Sec:Differentiability}$.
\section{Equivalence of the weak and strong formulation}\label{Sec:Weak} 
\subsection{Some known results on finite perimeter sets}
\begin{Def} We say that a sequence of finite perimeter sets $E_j$ \textbf{converges in $L^1_{loc}(M)$} to another finite perimeter set $E$, and we denote this by writing $E_j\rightarrow E$ in $L^1_{loc}(M)$, if $\chi_{E_j}\rightarrow\chi_{E}$ in $L^1_{loc}(M)$, i.e., if $V((E_j\Delta E)\cap U)\rightarrow 0\;\forall U\subset\subset M$. Here $\chi_{E}$ means the characteristic function of the set $E$ and the notation $U\subset\subset M$ means that $U\subseteq M$ is open and $\overline{U}$ (the topological closure of $U$) is compact in $M$.
\end{Def}
 \begin{Def}
We say that a sequence of finite perimeter sets $E_j$ \textbf{converge in the sense of finite perimeter sets} to another finite perimeter set $E$ if $E_j\rightarrow E$ in $L^1_{loc}(M)$, and 
\begin{eqnarray*} 
                           \lim_{j\rightarrow+\infty}\mathcal{P}(E_j)=\mathcal{P}(E).
\end{eqnarray*}
\end{Def}
For a more detailed discussion on locally finite perimeter sets and functions of bounded variation on a Riemannian manifold, one can consult \cite{MPPP}, for the more classical theory in $\R^n$ we refer the reader to \cite{AFP}, \cite{Maggi}.
\begin{Thm}[Fleming-Rishel]\label{Fle-Rish}
Let $u\in BV(M)$. Then the function $t\mapsto\mathcal{P}_{M}\left(\left\lbrace x\in M: u(x)>t \right\rbrace \right)$ is Lebesgue measurable on $\mathbb{R}$ and the following formula holds:
\begin{eqnarray}
\vert Du\vert(M)=\int_{-\infty}^{+\infty}\mathcal{P}_{M}\left(\left\lbrace x\in M: u(x)>t \right\rbrace \right)dt.
\end{eqnarray}
\end{Thm}
\begin{Dem}
See Theorem $4.3$ of \cite{ADo}.
\end{Dem}
\begin{Thm}[Proposition $1.4$ of \cite{MPPP}]\label{Prop:MPPPdensidade}
For every $u\in BV(M)$ there exists a sequence $(u_{j})_{j}\in C^{\infty}_{c}(M)$ such that $u_{j}\rightarrow u$ in $L^{1}_{loc}(M)$ and
\begin{eqnarray}
\vert Du\vert(M)=\lim_{j\rightarrow\infty}\int_{M}\vert\nabla u_{j}\vert dV_g.
\end{eqnarray}
\end{Thm}
\begin{Rem}
As a consequence of Theorem \ref{Prop:MPPPdensidade} we have
\begin{eqnarray}
\lim_{j\rightarrow\infty}\vert\left\lbrace x\in M: \vert u_{j}(x)-u(x)\vert\geq \eta\right\rbrace \vert=0,\,\,\,\,\, \forall\eta>0.
\end{eqnarray}
\end{Rem}
We state here a well known result.
\begin{Lemme}[Morse-Sard's Lemma]\label{suave} If $u\in C^{\infty}(M)$ and $E=\{x\in M: \nabla u(x)=0\}$, then $|u(E)|=0$. In particular, $\{u=t\}=\{x\in M: u(x)=t\}$ is a smooth hypersurface in $M$ for a.e. $t\in\mathbb{R}$.
\end{Lemme}
\subsection{Proof of the equivalence, Theorem $\ref{Thm:Equivalence}$}
Roughly speaking to prove Theorem $\ref{Thm:Equivalence}$ we make a construction which replace a finite perimeter set by one of the same volume with a small ball inside and one outside, by adding a small geodesic ball (with smooth boundary) to a point of density $0$ and subtracting a small geodesic ball to a point of of density $1$ taking care of not altering the volume. This enables us to obtain again a finite perimeter set of the same volume with a perimeter that is a small perturbation of the original one and that in addition have the property that we can put inside and outside a small ball. This construction legitimate us to apply mutatis mutandis the arguments of the proof of Lemma $1$ of \cite{Modica} to conclude the proof of Theorem $\ref{Thm:Equivalence}$. Our adapted version of Lemma $1$ of \cite{Modica} is the following lemma.
\begin{Lemme}\label{Lemma:1Modica} Let $\Omega\in\tilde{\tau}_M$, bounded, $\mathring{\Omega}\neq\emptyset$, and $Interior(\Omega^c)\neq\emptyset$. Then there exists a sequence $\Omega_k\in\tau_M$ with $V(\Omega_k)=V(\Omega)$, which converges to $\Omega$ in the sense of finite perimeter sets.
\end{Lemme}
\begin{Rem}
We observe that if $M$ is noncompact and $\Omega$ bounded, then we always have $Interior(\Omega^c)\neq\emptyset$. 
\end{Rem}
In connection with the original paper \cite{Modica}, we want just to point out two things. First, Lemma $1$ of \cite{Modica} is stated and proved in $\R^n$ but the proof generalizes immediately to complete Riemannian manifolds. The technical theorems needed to make this generalization are provided or are easily deducible from the paper \cite{MPPP} which extends the theory of $BV$-functions from $R^n$ to the setting of complete Riemannian manifolds. Second the assumption of $\Omega$ and $\Omega^c$ having nonvoid interior cannot be dropped to make the proof of Lemma $1$ of \cite{Modica} (and also Lemma \ref{Lemma:1Modica}) works. This is just a technical problem that we will solve in Lemma \ref{Lemma:Smoothisovolumic}.
The proof of Lemma $\ref{Lemma:1Modica}$ goes along the same lines of Lemma $1$ of \cite{Modica}, but to make the paper self contained we write it here. 

\begin{Dem}[ Lemma $\ref{Lemma:1Modica}$]
Take a bounded finite perimeter set $\Omega$ such that there exist $x_1\in\Omega$, $x_2\in M\setminus\Omega$, and $0<r_0<Min\{inj_{x_1}(M), inj_{x_2}(M)\}$ where for every $p\in M$, $inj_p(M)$ denotes the injectivity radius of $M$ at $p$, with $B_1:=B_M(x_1, r_0)\subseteq\Omega$ and $B_2:=B_M(x_2,r_0)\subseteq M\setminus\Omega$ and consider its characteristic function $\chi_{\Omega}$. Consider the usual mollifiers $\rho_{\varepsilon}$ described in Proposition $1.4$ of \cite{MPPP} and built the approximating functions $u_{\varepsilon}\in C^{\infty}_c(M)$ as described there.
By Proposition $1.4$ of \cite{MPPP} we have that $u_{\varepsilon}\to\chi_{\Omega}$ in $L^1(M)$ topology, when $\varepsilon\to 0^+$, in particular we have 
\begin{equation}\label{Eq:Convmeasure}
\lim_{\varepsilon\to0^+}V\left(\left\lbrace x\in M: |u_{\varepsilon}(x)-\chi_{E}(x)|\geq\eta\right\rbrace \right)=0, \forall\eta\geq 0. 
\end{equation} 
Moreover, $u_{\varepsilon}\in\mathcal{C}_{c}^{\infty}(M)$ and $\mathcal{P}_{M}(\Omega)=\lim_{\varepsilon\to 0^+}\int_{M}\vert\nabla u_{\varepsilon}\vert dV$. For every $\eta>0$ we can choose $0<\varepsilon_{\eta}<Min\{\eta,\frac{r_0}{2}\}$ such that 
$$V\left(\left\lbrace x\in M: |u_{\varepsilon_{\eta}}(x)-\chi_{E}(x)|\geq\eta\right\rbrace \right)\leq\eta.$$
Since $t\mapsto\mathcal{P}_{M}\left(\left\lbrace x\in M:u_{\varepsilon_{\eta}}(x)>t\right\rbrace  \right)$ is a Lebesgue measurable function, we can define 
$$A:=\{\mathcal{P}_{M}\left(\left\lbrace x\in M:u_{\varepsilon_{\eta}}(x)>t\right\rbrace  \right):\eta\leq t\leq 1-\eta\},$$ and 
\begin{equation}\label{Eq:Essinf}
\nu_{\eta}:=essinf_{\eta\leq t\leq 1-\eta}\mathcal{P}_{M}\left(\left\lbrace x\in M:u_{\varepsilon_{\eta}}(x)>t\right\rbrace  \right).
\end{equation}
By the very definition of $\nu_{\eta}$ we know that for every $\eta'>0$ we have $|[\nu_{\eta}, \nu_{\eta}+\eta']\cap A|>0$, thus using the Morse-Sard's theorem, i.e., Lemma $\ref{suave}$ and \eqref{Eq:Essinf} we get the existence of $t_{\eta}\in\left]\eta;1-\eta\right[$ such that
\begin{equation}\label{Eq:Essinf1}
\mathcal{P}_{M}\left( \left\lbrace x\in M:u_{\varepsilon_{\eta}}(x)>t_{\eta} \right\rbrace \right)<\nu_{\eta}+\eta,
\end{equation}
$t_{\eta}$ is a regular value of $u_{\varepsilon_{\eta}}$, i.e.,
$$\nabla u_{\varepsilon_{\eta}}(x)\neq 0,\:\forall x\in M: u_{\varepsilon_{\eta}}=t_{\eta}.$$
In view of this we can define $\Omega_{\varepsilon_{\eta}}'\in\tau_M$, $\Omega_{\varepsilon_{\eta}}':=u_{\varepsilon_{\eta}}^{-1}(]t_{\varepsilon_{\eta}},+\infty[)$, this ensures that $\Omega_{\varepsilon_{\eta}}'$ is bounded if $\Omega$ is bounded, furthermore we have also that $\partial\Omega_{\varepsilon}'=\left\lbrace x\in M: u_{\varepsilon_{\eta}}(x)=t_{\eta}\right\rbrace$ is smooth (again by Lemma $\ref{suave}$), and $\Omega_{\varepsilon_{\eta}}'\bigtriangleup\Omega\subset\left\lbrace x\in M: |u_{\varepsilon_{\eta}}(x)-\chi_{\Omega}(x)|\geq\eta\right\rbrace$. 
This last property joint with \eqref{Eq:Convmeasure} imply
\begin{eqnarray}\label{vol}
V(\Omega_{\varepsilon_{\eta}}'\bigtriangleup\Omega)\to 0,
\end{eqnarray}
which by the lower semicontinuity of the perimeter, gives 
\begin{eqnarray*}
\mathcal{P}_{M}(\Omega)\leqslant\liminf_{\eta\to0^+}\mathcal{P}_{M}(\Omega_{\varepsilon_{\eta}}').
\end{eqnarray*}
For the converse inequality, we deduce from \eqref{Eq:Essinf1} that
\begin{eqnarray*}
\mathcal{P}_{M}(\Omega'_{\varepsilon_{\eta}})\leqslant\eta+\nu_{\eta}\leqslant\eta+\mathcal{P}_{M}(\{x\in M: u_{\varepsilon_{\eta}}>t\}),
\end{eqnarray*}
for every $\eta>0$, and for almost all $t\in[\eta,1-\eta]$, so, integrating over the interval $[\eta,1-\eta]$ and applying the Fleming-Rishel formula (compare Theorem \ref{Fle-Rish}), we obtain 
\begin{eqnarray}
(1-2\eta)\mathcal{P}_{M}(\Omega_{\varepsilon_{\eta}})\leqslant\eta(1-2\eta)+\int_{M}|\nabla u_{\varepsilon_{\eta}}|dV,
\end{eqnarray}
which combined with Proposition \ref{Prop:MPPPdensidade} yields
\begin{eqnarray*}
\limsup_{\eta\to0^+}\mathcal{P}_{M}(\Omega_{\varepsilon_{\eta}})\leqslant\mathcal{P}_{M}(\Omega).
\end{eqnarray*}
Therefore we have proved that corresponding to every sequence $\eta_k\to0^+$, there exists a sequence $\Omega_k'\in\tau_M$ such that $\Omega_k'$ converges to $\Omega$ in the sense of finite perimeter sets. As it is easy to check for every $k$ large enough $B_1\subseteq\Omega_k'$ and $B_2\subseteq M\setminus\Omega_k'$. Set $\delta_k:=V(\Omega)-V(\Omega_k')$ and take $k$ large enough to ensure that $Min\{V(B_1),V(B_2)\}>|\delta_k|$.  
Now we choose $\Omega_{k}:=\Omega_k'\mathring{\cup}B(x_1, r_k)$, where $V(B(x_1, r_k))=|\delta_k|$,  if $\delta_k>0$, and $\Omega_{k}:=\Omega_k'\setminus B(x_2, r_k)$, where $V(B(x_2, r_k))=|\delta_k|$, if $\delta_k<0$, and finally $\Omega_k:=\Omega_k'$, if $\delta_k=0$. Using the fact that $\delta_k\to0$ we see that also $A(\partial B(x_i, r_k))\to 0$. It is straightforward to verify that $V(\Omega_k)=V(\Omega)$, $\partial\Omega_k$ is $C^{\infty}$, $\Omega_k$ is still bounded and 
$$V(\Omega_k\Delta\Omega'_k)\leq|\delta_k|\to0,\:k\to+\infty,$$
$$|\P(\Omega_k)-\P(\Omega_k')|\leq A(\partial B(x_i, r_k))\to0,\:k\to+\infty.$$ 
From the last properties it follows easily that the sequence $(\Omega_k)$ converges to $\Omega$ in the sense of finite perimeter sets, and the lemma follows.
\end{Dem}
We list here some lemmas that will be used in the proof of Lemma $\ref{Lemma:Smoothisovolumic}$.
\begin{Lemme}\label{cut}
Let $\Omega\subset M$ be any measurable set , then for all $J_{k}:=(k,2k+1)\subset (0,+\infty)$, $k\in\mathbb{N}$, there exists $r_{k}\in(k,2k+1)$ such that
$$\mathcal{H}^{n-1}(\Omega\cap\partial B_M(x,r_{k}))\leqslant \frac{V(\Omega)}{k},$$ where $x\in M$ is being taken fixed.
\end{Lemme}
\begin{Dem}
By coarea formula 
\begin{eqnarray*}
V(\Omega)=\int_{0}^{\infty}\mathcal{H}^{n-1}(\Omega\cap\partial B(x,r))dr,
\end{eqnarray*}
where $x$ is any fixed point in $M$.\\
We affirm that given $k\in\mathbb{N}$, there exists $r_{k}\in(k,2k+1)$ such that
\begin{eqnarray*}
\mathcal{H}^{n-1}(\Omega\cap\partial B(x,r_{k}))\leqslant \frac{V(\Omega)}{k},
\end{eqnarray*}
because otherwise we would have 
\begin{eqnarray*}
V(\Omega)\geqslant\int_{k}^{2k+1}\mathcal{H}^{n-1}(\Omega\cap\partial B(x,r))dr>\frac{(k+1)V(\Omega)}{k},
\end{eqnarray*}
which is a contradiction.
\end{Dem}
\begin{Rem}
See that when $r_{k}\to\infty$, it holds
\begin{eqnarray*}
V(\Omega\cap B(x,r_{k}))\to V(\Omega),\,\,\,\, k\to\infty.
\end{eqnarray*}
\end{Rem}
\begin{Lemme}\label{Lemma:Smoothisovolumic} Let $\Omega\in\tilde{\tau}_M$ be a finite perimeter set with $V(\Omega)<+\infty$, $V(\Omega), V(\Omega^c)>0$, where $\Omega^c:=M\setminus\Omega$. Then there exists a sequence $\Omega_{k}\in\tau_{M}$ such that $V(\Omega_k)=V(\Omega)$ and $\Omega_k$ converges to $\Omega$ in the sense of finite perimeter sets. 
\end{Lemme}
In the proof of this Lemma we really differ from the paper \cite{Modica}, even if we make a crucial use of Lemma $1$ of that paper.
 
\begin{Dem}
Consider an arbitrary set $\Omega\in\tilde{\tau}_M$ and take two distinct points $x_1\in\Omega$ and $x_2\in\Omega^c$ of density $\Theta(x_1, V\llcorner\Omega)=1$ and $\Theta(x_2, V\llcorner\Omega)=0$, where $\Theta(p, V\llcorner\Omega):=\lim_{r\to0^+}\frac{V(\Omega\cap B(p, r))}{\omega_nr^n}$, for every $p\in M$.  By $\omega_n$ we denote the volume of the ball of radius $1$ in $\R^n$. Consider the two continuous functions $f_1, f_2:I\to\R$, where $I:=[0, r_0[$ such that $f_1(r):=V(\Omega\cap B_M(x_1, r))$, $f_2(r):=V(\Omega^c\cap  B_M(x_2, r))$. The radius $r_0$ could be chosen small enough to have $B_M(x_1, r_1)\cap B_M(x_2, r_2)=\emptyset$ for every $r_1, r_2\in I$ and such that there exist $r_1, r_2\in I$ satisfying the property $f_1(r_1)=f_2(r_2)$ and $\partial B_M(x_1, r_1),\partial B_M(x_2, r_2)$ smooths (for this last property it is enough to take $r_0$ less than the injectivity radius at $x_1$ and $x_2$). Then we set 
$$\tilde{\Omega}:=[\Omega\setminus B_M(x_1,r_1)]\mathring{\cup}[\Omega^c\cap B_M(x_2, r_2)]=[\Omega\setminus B_M(x_1,r_1)]\cup B_M(x_2, r_2).$$ 
As it is easy to see $V(\tilde{\Omega})=V(\Omega)$, 
\begin{equation}\label{Eq:Equivalence}
|\P(\tilde{\Omega})-\P(\Omega)|\leq \sum_{i=1}^2[A(\partial B_M(x_i, r_i))+\P(\Omega, B_M(x_i, r_i))],
\end{equation}
\begin{equation}\label{Eq:Equivalence1}
V(\Omega\Delta\tilde{\Omega})=f_1(r_1)+f_2(r_2),
\end{equation} $\mathring{\tilde{\Omega}}\neq\emptyset$, and $Interior(\tilde{\Omega}^c)\neq\emptyset$. It is straightforward to verify that the right hand sides of \eqref{Eq:Equivalence} and \eqref{Eq:Equivalence1} converges to zero when the radii $r_1$ and $r_2$ go to zero. We prove the lemma first for bounded sets $\Omega\in\tilde{\tau}_M$, and then we pass to the general case by observing that one can always approximate an unbounded $\Omega\in\tilde{\tau}_M$ in the sense of finite perimeter sets by a sequence of bounded ones. Let us assume that $\Omega\in\tilde{\tau}_M$ is bounded, then for any arbitrary $\varepsilon>0$, the Riemannian version of Lemma $1$ of \cite{Modica}, namely Lemma  $\ref{Lemma:1Modica}$ applied to $\tilde{\Omega}$ permits to find $\tilde{\Omega}_{\varepsilon}\in\tau_M$ such that $V(\tilde{\Omega}_{\varepsilon})=V(\tilde{\Omega})=V(\Omega)$ and $$V(\tilde{\Omega}_{\varepsilon}\Delta\tilde{\Omega})\leq\frac{\varepsilon}{2},$$ 
$$|\P(\tilde{\Omega}_{\varepsilon})-\P(\tilde{\Omega})|\leq\frac{\varepsilon}{2}.$$
These last two inequalities combined with \eqref{Eq:Equivalence} and \eqref{Eq:Equivalence1} imply that
\begin{equation}
V(\tilde{\Omega}_{\varepsilon}\Delta\Omega)\leq\varepsilon,
\end{equation} 
\begin{equation} 
|\P(\tilde{\Omega}_{\varepsilon})-\P(\Omega)|\leq\varepsilon.
\end{equation}
To finish the proof we consider now the case of an unbounded $\Omega\in\tilde{\tau}_M$, with $V(\Omega)=v<+\infty$. Fix a point $p\in M$, a fine use of the coarea formula as explained in Lemma \ref{cut} gives a sequence of radii $r_{k}\to\infty$, $r_{k}\geqslant k$, such that whenever $B(p,r_{k})\cap\Omega=:\Omega_{k}$ we have
 \begin{eqnarray*}
 \lim_{k\to\infty}\mathcal{P}(\Omega_{k})=\mathcal{P}(\Omega),
 \end{eqnarray*}
 because
 \begin{eqnarray*}
\mathcal{P}(\Omega_{k})\leqslant\mathcal{P}(\Omega,B(p,r_{k}))+\frac{V(\Omega)}{k}, 
 \end{eqnarray*}
 which after taking limits leads to $$\limsup\mathcal{P}(\Omega_{k})\leqslant\mathcal{P}(\Omega),$$
 and because from $V(\Omega\bigtriangleup\Omega_{k})\to 0$ and the lower semicontinuity of the perimeter we get $\liminf\mathcal{P}(\Omega_{k})\geqslant\mathcal{P}(\Omega)$. Now, we observe that $V(\Omega_k)\leq V(\Omega)$ and in general it could happen that $V(\Omega_k)<V(\Omega)$, so  it still remains to readjust the volumes of these $\Omega_k$'s. We do it by perturbing $\Omega_k$ in adding a small geodesic ball $B_M(p_1,r'_k)$ such that $V(B_M(p_1, r'_k)\cap\Omega^c)=v-v_k$, with $v_{k}=V(\Omega_{k})$, centered at a fixed point $p_1$ of density $\Theta(p_1, V\llcorner\Omega)=0$, with $k$ sufficiently large. It is worth to note that as above (see Lemma $\ref{Lemma:1Modica}$) $r'_k\to0$, when $k\to\infty$. This construction gives sets $\Omega_k'\in\tilde{\tau}_M$, such that $V(\Omega_{k}')=v=V(\Omega)$, $\Omega_k'$ is bounded, 
 \begin{equation}\label{Eq:Equivalence3}
 V(\Omega_k\Delta\Omega_k')=v-v_k,
 \end{equation} 
 \begin{equation}\label{Eq:Equivalence4}
 |\mathcal{P}(\Omega_{k}^{'})-\mathcal{P}(\Omega_{k})|\leq\P(\Omega, B_M(p_1,r'_k))+A(\partial B_M(p_1,r'_k)).
 \end{equation}
 Therefore $(\Omega_k')_{k\in\mathbb{N}}$ converges to $\Omega$ in the sense of sets of finite perimeter, because the right hand sides of \eqref{Eq:Equivalence3} and \eqref{Eq:Equivalence4} go to $0$, when $k\to+\infty$. Unfortunately the sets $\Omega_k'$ are still not open with smooth boundary. Hence we have to continue our construction to achieve the proof. To do so consider any given sequence of positive numbers $\varepsilon_k\to0$ the fact that $\Omega_k'$ is bounded allow us, as above in Lemma $\ref{Lemma:Smoothisovolumic}$, to find $\Omega_k''\in\tau_M$ such that  $V(\Omega_k'')=V(\Omega_k')=V(\Omega)=v$
 \begin{equation}
 V(\Omega_k'\Delta\Omega_k'')\leq\varepsilon_k,
 \end{equation}  
 \begin{equation} 
 |\P(\Omega_k')-\P(\Omega_k'')|\leq\varepsilon_k.
 \end{equation}
Since the sequence $(\Omega_k')_{k\in\mathbb{N}}$ converges to $\Omega$ in the sense of sets of finite perimeter, the last two equations ensures that the sequence $\Omega_k''\in\tau_M$ converges to $\Omega$ in the sense of finite perimeter sets too, which is our claim. 
\end{Dem}
Now we are ready to prove easily Theorem $\ref{Thm:Equivalence}$.

\begin{Dem}[of Theorem $\ref{Thm:Equivalence}$]
Taking into account Remark $\ref{Rem:Trivialinequality}$, to show the theorem, it is enough to prove the nontrivial inequality $I_M(v)\leq\tilde{I}_M(v)$ for every $v\in[0, V(M)[$.  To this aim, let us consider $\varepsilon>0$ and $\Omega\in\tilde{\tau}_M$, with $V(\Omega)=v$. By Lemma $\ref{Lemma:Smoothisovolumic}$ there is a sequence $\Omega_k\in\tau_M$ such that $V(\Omega_k)=v$, and $(\Omega_k)$ converging to $\Omega$ in the sense of finite perimeter sets. In particular we have that $\lim_{k\to+\infty}\P(\Omega_k)=\P(\Omega)$. On the other hand by definition we have that 
$I_M(v)\leq\P(\Omega_k)$ for every $k\in\mathbb{N}$. Passing to limits leads to have 
\begin{equation}\label{Eq:Thm1}
 I_M(v)\leq\P(\Omega),
\end{equation} 
for every $\Omega\in\tilde{\tau}_M$ with $V(\Omega)=v$. Taking the infimum in \eqref{Eq:Thm1} when $\Omega$ runs over $\tilde{\tau}_M$ keeping $V(\Omega)$ fixed and equal to $v$, we infer that $I_M(v)\leq\tilde{I}_M(v)$. This completes the proof.
\end{Dem}

\section{Continuity of $I_M$}\label{Sec:Continuity}  
\subsection{Continuity in bounded geometry}\label{Subsec:Continuity}
To illustrate the proof of Theorem \ref{Main} we start this section with the easy part of the proof resumed in the next lemma that is straightforward, compare \cite{MorganBlog} Proposition $1$. As the example $3.53$ of \cite{AFP} shows, in general we can have finite perimeter sets with positive perimeter and void interior that are not equivalent to any other set of finite perimeter with non void interior. So the question of putting a ball inside or outside a set of finite perimeter is a genuine technical problem. On the other hand, following \cite{Tamanini} Theorem $1$, it is always possible to put a small ball inside and outside an isoperimetric region, which justify the constructions performed in this proof. As a general remark a result of Federer (the reader could consult \cite{AFP} Theorem $3.61$) states that for a given set of finite perimeter $E$ the density is either $0$ or $\frac{1}{2}$ or $1$, $\mathcal{H}^{n-1}$-a.e. $x\in M$, moreover points of density $1$ always exist $V$-a.e. inside $D$, because of the Lebesgue's points Theorem applied to the characteristic function of any $V$-measurable set of $M$. About this topic the reader could consult the book \cite{Maggi} Example $5.17$. Thus $V(D)>0$ ensures the existence of at least one point $p$ belonging to $D$ of density $1$, which is enough for the aims of our proofs.
\begin{figure}[h!]
  \centering
    \includegraphics[width=.5\linewidth]{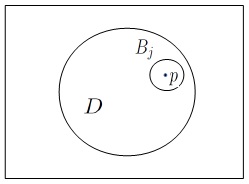}
  \caption{$v'<v$ Upper Semicontinuity}
  \label{fig:fig1}
\end{figure}
\begin{figure}[h!]
  \centering
    \includegraphics[width=.5\linewidth]{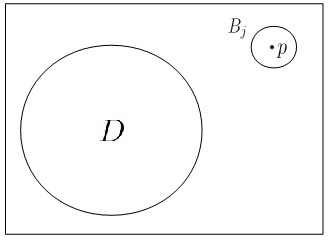}
  \caption{$v'>v$ Upper Semicontinuity}
  \label{fig:fig2}
\end{figure}
\begin{Thm}\label{Thm:Uppersemicontinuity} Let $M$ be a Riemannian manifold (possibly incomplete, or possibly complete not necessarily with bounded geometry). If there exists an isoperimetric region in volume $v\in]0, V(M)[$, then $I_M$ is upper semicontinuous in $v$. 
\end{Thm}
\begin{Dem} To prove the theorem it is enough to prove the next two inequalities.
\begin{equation}\label{Eq:Cont3bis}
\limsup_{v'\rightarrow v^-} I_M(v')\leq I_M(v).
\end{equation}
 \begin{equation}\label{Eq:Cont4bis}
\limsup_{v'\rightarrow v^+} I_M(v')\leq I_M(v).
\end{equation}
In first we prove (\ref{Eq:Cont3bis}). If  $v_j\nearrow v$, consider an isoperimetric region $D$ in volume $V(D)=v$, 
$$I_M(v)=A(\partial D).$$ 

Then for $j$ sufficiently large one can subtract a small geodesic ball (i.e. of small radius) $B_j=B(p, r'_j)$ of volume $v-v_j$ from $D$, centered at some point $p\in D$ of density $1$, to obtain $D'_j:=D\setminus B(p,r'_j)$ of volume $V(D'_j)=v_j$ and $A(\partial D'_j)\leq A(\partial D)+A(\partial B_j)$. Observe here that the center $p$ of $B_j$  is fixed with respect to $j$. Moreover $r'_j\rightarrow 0$, and this is always possible to obtain in any Riemannian manifold. So by definition of $I_M(v_j)$,  holds $$I_M(v_j)\leq A(\partial D'_j)\leq A(\partial D)+A(\partial B_j)=I_M(v)+A(\partial B_j),$$ which implies that $$\limsup I_M(v_j)\leq\limsup A(\partial D)+A(\partial B_j)\leq I_M(v),$$
since the sequence $v_j$ is arbitrary we get (\ref{Eq:Cont3bis}). In second, we prove (\ref{Eq:Cont4bis}). 
If $v_j\searrow v$, then take an isoperimetric region of voume $v$, i.e., $V(D)=v$, $A(\partial D)=I_M(v)$ and then add a small ball $B_j:=B(p, r_j)$ of volume $v_j-v$ to $D$ outside $D$ to obtain $D'_j:=D\mathring{\cup}B_j$ of volume $V(D'_j)=v_j$ and $A(\partial D'_j)=A(\partial D)+A(\partial B_j)$. Observe again that the center $p$ of $B_j$ here is fixed with respect to $j$ and $r_j\rightarrow 0$, this is always possible in any Riemannian manifold. By definition of $I_M(v_j)$ we get $$I_M(v_j)\leq A(\partial D'_j)=A(\partial D)+A(\partial B_j)=I_M(v)+A(\partial B_j),$$ now taking the $\limsup$ it follows $$\limsup I_M(v_j)\leq\limsup [A(\partial D)+A(\partial B_j)]=I_M(v)+\limsup A(\partial B_j)=I_M(v),$$
since the sequence $v_j$ is arbitrary we get (\ref{Eq:Cont4bis}), which completes the proof.
\end{Dem}

At this point, we may finish the proof of the main Theorem $\ref{Main}$. 
\begin{Dem}[of Theorem $\ref{Main}$] We will prove separately the following four inequalities that together will give the proof of our Theorem $\ref{Main}$. 
\begin{equation}\label{Eq:Cont1}
I_M(v)\leq\liminf_{v'\rightarrow v^-} I_M(v').
\end{equation}
\begin{equation}\label{Eq:Cont2}
I_M(v)\leq\liminf_{v'\rightarrow v^+}I_M(v').
\end{equation}
\begin{equation}\label{Eq:Cont3}
\limsup_{v'\rightarrow v^-} I_M(v')\leq I_M(v).
\end{equation}
 \begin{equation}\label{Eq:Cont4}
\limsup_{v'\rightarrow v^+} I_M(v')\leq I_M(v).
\end{equation}
\begin{figure}[h!]
  \centering
    \includegraphics[width=.5\linewidth]{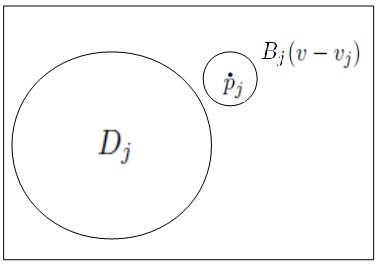}
  \caption{$v'<v$ Lower Semicontinuity}
  \label{fig:fig1}
\end{figure}

To prove (\ref{Eq:Cont1}) we want to add a small ball. Let $v_j\nearrow v$, take a domain $D_j$ in volume $v_j$ such that $V(D_j)=v_j$ and $I_M(v_j)\leq A(\partial D_j)+\frac{1}{j}$, then add a small ball $B_j:=B(p_j, r_j)$ to $D_j$ outside $D_j$ to obtain $D'_j$ of volume $v$ and $A(\partial D'_j)=A(\partial D_j)+A(\partial B_j)$.  This is possible because $D_j$ by the very definition (see Definition \ref{Def:IsPStrong}) may be chosen bounded. It is worth to observe here that the centers $p_j$ are variable and not fixed as in the proof of Theorem \ref{Thm:Uppersemicontinuity}. So we need to use Bishop-Gromov's Theorem to bound the area of $B_j$ uniformly w.r.t. the centers. Having in mind the definition of $I_M(v)$ it is easy to see that
\begin{equation} I_M(v)=I_M(V(D'_j))\leq A(\partial D_j)=A(\partial D_j)+A(\partial B_j).
\end{equation} 
Now observe that by Lemma $3.2$ of \cite{MonNar} or Lemma $3.5$ of \cite{MJ} that $A(\partial B_j)\leq A(\partial B_{\mathbb{M}^n_k}(v-v_j))$ where 
$B_{\mathbb{M}^n_k}(w)$ is a geodesic ball enclosing volume $w$ in $\mathbb{M}^n_k$. As it is easy to check $A(\partial B_{\mathbb{M}^n_k}(w))\rightarrow 0$, when $w\rightarrow 0$, because the centers could be chosen fixed in the comparison manifold. This implies that $A(\partial B_{\mathbb{M}^n_k}(v-v_j))\rightarrow 0$, when $j\rightarrow+\infty$ and a fortiori that $\liminf_{j\rightarrow+\infty} A(\partial B_j)=0$. Thus 
\begin{equation} I_M(v)\leq A(\partial D'_j)\leq I_M(v_j)+\frac{1}{j}+A(\partial B_{\mathbb{M}^n_k}(v-v_j))\leq\liminf\: I_M(v_j).
\end{equation} 
By the arbitrariness of the initial sequence of volumes $(v_j)$, (\ref{Eq:Cont1}) follows readily.

To show (\ref{Eq:Cont2}) the strategy is now to subtract a small ball to an eventually diverging (to infinity) sequence of domains that could become thinner and thinner without leaving the opportunity of placing a small ball of the right value of the volume inside them. To rule out  this possibility Lemma 2.5 of \cite{Nar12} is needed. This is a more delicate task with respect to the preceding construction in which we add a small ball to a relatively compact domain. 
\begin{figure}[h!]
  \centering
    \includegraphics[width=.8\linewidth]{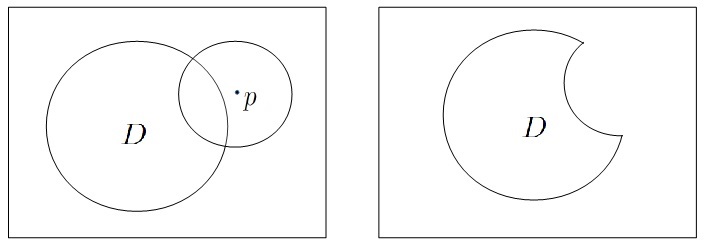}
  \caption{$v'>v$ Lower Semicontinuity}
  \label{fig:fig1}
\end{figure}
\begin{Rem} From the proof of Lemma 2.5 of $\cite{Nar12}$ we argue that when $|v-v'|\sim r^n<<v$, $m'_0=\frac{1}{2}c_1(n,k,r)=\frac{r^n}{2e^{(n-1)\sqrt{|k|}}}$.\end{Rem} 

Let $D\in\tau_M$, (this means that $D$ is open, $D\subset\subset M$, and $\partial D$ is $C^{\infty}$) such that $V(D)=v'>v$ and then take $r$ satisfying $\frac{r^nv_0}{2e^{(n-1)\sqrt{k}}}=v'-v$ (this is possible because the function $r\mapsto V(B(p,r)\cap D)$ is continuous whenever $D$ is a $V$-measurable set, as it is easy to check), by Lemma 2.5 of \cite{Nar12} we may take a point $p\in M$ such that for small $v'-v$ one have
\begin{equation}\label{Eq:Main0}
     V(B(p,r)\cap D)>\frac{r^nv_0}{2e^{(n-1)\sqrt{k}}}=v'-v.
\end{equation} 
This is possible because for small $|v-v'|$ we can take $r$ small enough to obtain that the constant $m'_0$ produced by Lemma 2.5 of \cite{Nar12} coincides with the right hand side of the preceding inequality. An easy consequence of  (\ref{Eq:Main0}) is that
$$V(D\setminus B(p,r))=V(D)-V(B(p,r)\cap D)<v,$$
it follows that we may choose $0<r'< r$ satisfying $V(D\setminus B(p,r'))=v$, because the function $r\mapsto V(D\setminus B(p, r)$ is continuous.  Fix $\eta>0$ and consider an almost isoperimetric region $D\in\tau_M$ in volume $v'$, i.e., such that $V(D)=v'$ and
\begin{equation}
I_M(v')\leq A(\partial D)\leq I_M(v')+\eta,
\end{equation}  
by Bishop-Gromov's theorem it is true that $A(\partial B_M(p, r'))\leq A(\partial B_{\mathbb{M}^n_k}(r'))$, then by Theorem $\ref{Thm:Equivalence}$ we have the following
\begin{eqnarray}\label{Eq:Main1}
I_M(v) & \leq & A(\partial (D\setminus B_M(p, r')))\leq A(\partial D)+A(\partial B_M(p, r'))\\ \label{Eq:Main1bis}
& \leq & I_M(v')+\eta+A(\partial B_{\mathbb{M}^n_k}(r')), 
\end{eqnarray}
with $r'<r=\left(2\frac{v'-v}{v_0}e^{(n-1)\sqrt{k}}\right)^{\frac{1}{n}}$. 
\begin{Rem}
 The second inequality in $(\ref{Eq:Main1})$ holds also more generally for sets of locally finite perimeter. For a rigorous discussion of this fact the reader could consult Proposition $3.38$ $(d)$ of $\cite{AFP}$, Lemma $12.2$, and Exercise $12.23$ of $\cite{Maggi}$. 
\end{Rem} 
\begin{Rem}\label{Rem:Crucial} It is worth to note that the competitor $D\setminus B(p, r')$ have a not necessarily smooth boundary, but still $D\setminus B(p, r')\in\tilde{\tau}_M$. We need Theorem $\ref{Thm:Equivalence}$ to completely justify the first inequality of \eqref{Eq:Main1}. Alternatively we can observe that when $V(D\setminus B(p, r'))>0$, $D\setminus B(p, r')$ and its complement have nonvoid interior and we can approximate it simply by Lemma $\ref{Lemma:1Modica}$ without using the entire strenght of Theorem $\ref{Thm:Equivalence}$. This latter argument still justify \eqref{Eq:Main1}. 
\end{Rem}
By the arbitrariness of $\eta>0$ we get
\begin{equation}
I_M(v)\leq I_M(v')+A(\partial B_{\mathbb{M}^n_k}(r')).
\end{equation}
Taking limits in the last inequality yields 
\begin{equation}
I_M(v)\leq\liminf_{v'\rightarrow v^+}I_M(v').
\end{equation}

The last two inequalities are relative to the $\limsup$ property and are analogous to the case in which there is existence of an isoperimetric region of voume $v$, but with the additional difficulty that isoperimetric regions of volume $v$ does not necessarily exists. So we apply the same ideas of the proof of Theorem \ref{Thm:Uppersemicontinuity} to a minimizing sequence of volume $v$ instead of a genuine isoperimetric region. 

Now, we prove (\ref{Eq:Cont3}). If  $v_j\nearrow v$, consider an almost minimizer $D_j\in\tau_M$ of volumes $V(D_j)=v$, i.e., 
$$I_M(v)\leq A(\partial D_j)\leq I_M(v)+\frac{1}{j}.$$ 
\begin{figure}[h!]
  \centering
    \includegraphics[width=.8\linewidth]{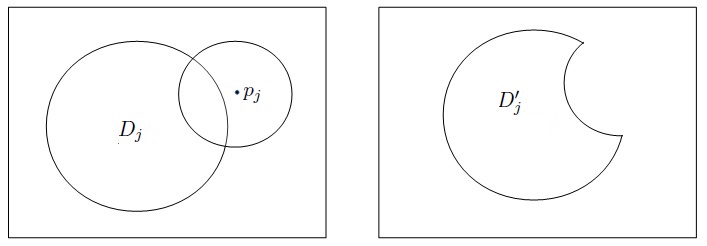}
  \caption{$v'<v$ Upper Semicontinuity}
  \label{fig:fig1}
\end{figure}
Then subtract a small ball $B_j:=B(p_j,r'_j)$ (whose intersection with $D_j$, $B_j\cap D_j$ has volume $v-v_j$) to $D_j$ as in the proof of (\ref{Eq:Cont2}), to obtain $D'_j:=D_j\setminus B(p_j,r'_j)$ of volume $V(D'_j)=v_j<v$ and $$A(\partial D'_j)\leq A(\partial D_j)+A(\partial B_j),$$ so by definition and Theorem  $\ref{Thm:Equivalence}$, (see Remark $\ref{Rem:Crucial}$) it holds $$I_M(v_j)\leq A(\partial D'_j)\leq A(\partial D_j)+A(\partial B_j),$$ which implies (as in the proof of (\ref{Eq:Cont2})) that $$\limsup I_M(v_j)\leq\limsup [A(\partial D_j)+A(\partial B_j)]= I_M(v).$$
Since the sequence $(v_j)$ is arbitrary we get (\ref{Eq:Cont3}).

Finally we prove (\ref{Eq:Cont4}). This last part of the proof is analogous in some respects to the proof of (\ref{Eq:Cont1}), because we add a small ball. 
\begin{figure}[h!]
  \centering
    \includegraphics[width=.5\linewidth]{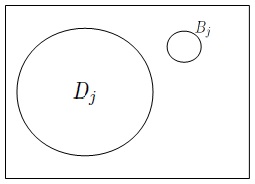}
  \caption{$v'>v$ Upper Semicontinuity}
  \label{fig:fig1}
\end{figure}
If $v_j\searrow v$, then take a minimizing sequence $D_j$ of volume $v$, i.e., $V(D_j)=v$, $A(\partial D_j)\searrow I_M(v)$ and then add a small ball $B_j$ to $D_j$ outside $D_j$ to obtain $D'_j$ of volume $V(D'_j)=v_j$ and $A(\partial D'_j)=A(\partial D_j)+A(\partial B_j)$, $$I_M(v_j)\leq A(\partial D'_j)=A(\partial D_j)+A(\partial B_j),$$ now taking the $\limsup$ it follows as before that $$\limsup \ I_M(v_j)\leq\limsup \ A(\partial D_j)+A(\partial B_j)=I_M(v)+\limsup \ A(\partial B_j)=I_M(v),$$
since the sequence $(v_j)$ is arbitrary we get (\ref{Eq:Cont4}), which completes the proof.
\end{Dem}
\section{Differentiability of $I_M$}\label{Sec:Differentiability}
\begin{Lemme}[Lemma $3.2$ of \cite{MJ}]\label{Lem:BVMJ1} Let $f:]a, b[\rightarrow\R$ be a continuous function. Then $f$ is concave (resp. convex) if and only if for every $x_0\in ]a, b[$ there exists an open interval $I_{x_0}\subseteq ]a,b[$ of $x_0$ and a concave (resp. convex) $C^2$ function $g_{x_0}:I_{x_0}\rightarrow\R$ such that $g_{x_0}=f(x_0)$ and $f(x)\leq g_{x_0}(x)$ (resp. $f(x)\geq g_{x_0}(x)$) for every $x\in I_{x_0}$.
\end{Lemme}
\begin{Rem} The preceding Lemma is just a rephrasing of the supporting hyperplanes theorem for closed convex sets of $\R^n$. To apply it in our context the hypothesis of continuity is crucial, we cannot assume $f$ just upper (resp. lower) semicontinuous. In fact take as a counterexample a function that is strictly monotone increasing on $[a, b]$, right continuous in an interior point $x_0$ but not continuous at $x_0$ with a strictly positive jump in $x_0$, concave at the left of $x_0$ and to the right of $x_0$. This function is not concave on the entire interval $[a, b]$, is upper semicontinuous and satisfies the other hypothesis of Lemma \ref{Lem:BVMJ1}, except continuity.  
\end{Rem}
\begin{figure}[h!]
 \centering
    \includegraphics[width=.5\linewidth]{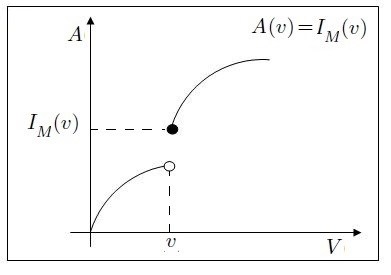}
  \caption{An example of an upper semicontinuous function that satisfies all the assumptions of Lemma \ref{Lem:BVMJ1} but the continuity.}
  \label{fig:fig1.4}
\end{figure}
We recall here the generalized existence Theorem $1$ of \cite{Nar12} stated under more general assumptions to check why this is legitimate one can see Remark 2.9 of \cite{MonNar}, or Remarks \ref{Rem:1}, \ref{Rem:LBRic}. 
\begin{Thm}[Generalized existence]\label{Thm:GenExistence}
Let $M$ have $C^0$-locally asymptotically bounded geometry in the sense of Definition \ref{Def:BoundedGeometryInfinity}. Given a positive volume $0<v < V(M)$, there are a finite number of limit manifolds at infinity such that their disjoint union with M contains an isoperimetric region of volume $v$ and perimeter $I_M(v)$. Moreover, the number of limit manifolds is at worst linear in $v$.
\end{Thm}
\begin{Rem}\label{Rem:1} The regularity discussion made there in Remark $2.2$ of \cite{MonNar}, is necessary in the proof of Corollary \ref{BPGen}, where we need to do analysis on the limit manifolds, applying a (by now classical) formula for the second variation of the area functional preserving volumes on those isoperimetric regions which eventually lie in a limit manifold of possibly non-smooth boundary.  The assumption of $C^0$ convergence of the metric tensor in the preceding lemma is due to the necessity of transporting volumes and perimeters in the limit manifold.
\end{Rem}
\begin{Rem}\label{Rem:LBRic}
We observe that if $(M_i,g_i,p_i)\rightarrow (M,g,p)$ in the pointed Gromov-Hausdorff topology and $M_i$ satisfy $Ric_{g_i}\geq (n-1) k_0 g_i$, it is not true, in general, that $Ric_{g}\geq (n-1) k_0 g$. Instead, if $(M_i,g_i,p_i)\rightarrow (M,g,p)$ in the pointed $C^0$-topology then $(M_i,g_i,V_i, p_i)\rightarrow (M,g,V, p)$ converge in the measured pointed Gromov-Hausdorff topology. Therefore, if all the Riemannian $n$-manifolds $(M_i,g_i)$ satisfy $Ric_{g_i}\geq (n-1) k_0 g_i$ then also the limit Riemannian manifold $(M,g)$  satisfies $Ric_g \geq (n-1) k_0 g$ (see Section 7 in $\cite{AG}$). Notice that for the convergence of the  Ricci curvature one should need a stronger convergence of the $(M_i,g_i,p_i)$ to $(M,g,p)$, say in $C^2$-topology; here we just need the convergence of a lower bound.
\end{Rem}
\begin{Rem} One possible application of Theorem $\ref{Thm:GenExistence}$, is to simplify part of the proof of different papers appeared in the literature about existence and characterisation of isoperimetric regions in noncompact Riemannian manifolds and prove new theorems of the same kind, as for instance it is done in Theorem $1$ of \cite{FloresNar016}.\end{Rem}
We can finish now the proof of Corollary \ref{BPGen}.
  
\begin{Dem} Using the generalized existence theorem of \cite{Nar12} and evaluating the second variation formula for the area functional on a generalized isoperimetric region $\Omega_{\bar{v}}$ in volume $V(\Omega_{\bar{v}})=\bar{v}$ we can construct a smooth function $f_{\bar{v}}$ defined in a small neighborhood of $\bar{v}$, that we can compare locally with $I_M$. Consider the equidistant domains $\Omega_t:=\left\{x\in M:\: d(x, \Omega_{\bar{v}})\leq t\right\}$, if $r_{\bar{v}}\geq t\geq 0$, and $\Omega_t:=M\setminus\left\{x\in M:\:d(x,M\setminus\Omega_{\bar{v}})\leq t\right\}$, if $-r_{\bar{v}}\leq t<0$, where $r_{\bar{v}}>0$ is the normal injectivity radius of $\partial\Omega_{\bar{v}}$. Consider the inverse function of $t\mapsto V(\Omega_t)$ as a function of the volume, $v\mapsto t(v)$, and finally set $f_{\bar{v}}(v):=A(\partial\Omega_{t(v)})$ for $v$ belonging to a small neighbourhood $I_{\bar{v}}=[\bar{v}-\varepsilon_{\bar{v}}, \bar{v}+\varepsilon_{\bar{v}}]$. To be rigorous in this construction we have to take care of the singular part of the boundaries of domains $\Omega_t$. This is done, carefully, in Proposition 2.1 and 2.3 of \cite{Bayle}. Here we just ignore this technical complication, to make the exposition simpler to read. We just observe that the proof that we give here works mutatis mutandis also if we consider the case in which $\Omega$ is allowed to have a nonvoid singular part. Hence, for every 
$\bar{v}\in]0, V(M)[$, $f_{\bar{v}}$ gives smooth function 
$f_{\bar{v}}:[\bar{v}-\varepsilon_{\bar{v}}, \bar{v}+\varepsilon_{\bar{v}}]\rightarrow[0, +\infty[$, such that $f_{\bar{v}}(\bar{v})=I_M(\bar{v})$ and $f_{\bar{v}}\geq I_M$. A standard application of the second variation formula  see (V.4.3) \cite{Chavel}, or \cite{BP}, shows that
\begin{equation}\label{Eq:BP}
f_{\bar{v}}''(v)=-\frac{1}{f_{\bar{v}}^2(v)}\left\{\int_{\partial\Omega_{t(v)}}(|II|^2+Ricci(\nu))d\H^{n-1}\right\}.
\end{equation}
From an elementary fact of linear algebra we know that $|II|^2\geq \frac{h^2}{n-1}$. Hence substituting in the preceding inequality, we get
\begin{equation} 
          f_{\bar{v}}''(v)\leq-\frac{(n-1)k}{f_{\bar{v}}(v)}.
\end{equation} 
If $k\geq 0$, then $f_{\bar{v}}$ is concave and a straightforward application of Lemma \ref{Lem:BVMJ1} implies that $I_M$ is concave in all $]0, V(M)[$.
If $k<0$ then
\begin{equation} 
          f_{\bar{v}}''(v)\leq-\frac{(n-1)k}{I_M(v)},
\end{equation}          
\begin{equation}
C=C(n,k, a,b):=\frac{(n-1)k}{2\delta_{M, a, b}},
\end{equation} 
where $\delta_{M, a, b}:=\inf\{I_M(v):v\in[a, b]\}$ is strictly positive because by Theorem \ref{Main}, $I_M$ is continuous. For every $\bar{v}\in ]a, b[$ it is easily seen that $$I_M(v)+C(a,b)v^2\leq f_{\bar{v}}(v)+C(a,b)v^2,$$ with $${(f_{\bar{v}}(v)+C(a,b)v^2)}''\leq 0,$$ for every $v\in ]a,b[\cap I_{\bar{v}}$. By Lemma \ref{Lem:BVMJ1}, for $a, b\in]0, V(M)[$, $I_M(v)+C(a,b)v^2$ is concave in $[a, b]$. Hence, $I_M(v)+C(a,b)v^2$ is locally Lipschitz and it is straightforward to see that $I_M$ is locally Lipschitz too, with ${I'}^+\leq {f'}_{\bar{v}}\leq{I'}^-$, with equality holding at all but a countable set of points, which are the only points of discontinuity of ${I'}^+$ and ${I'}^-$. Moreover ${I'}^+$ and ${I'}^-$ are nonincreasing so the set of points at which $I_M$ is nonderivable is at most countable, moreover $I'_M$ or $I'_M+2Cv$ are respectively monotone nonincreasing see for this standard convexity arguments Corollary $2$, page 29 of \cite{Bourbaki} this implies that they are special cases of absolutely continuous functions and for this reason differentiable almost everywhere. So exists $I''_M(v)$ almost everywhere. Now, following \cite{Bayle}, for an arbitrary function $f$, set 
\begin{equation}\label{Eq:Bayle}
\overline{D^2f}(x_0):=\limsup_{\delta\rightarrow0}\frac{f(x_0+\delta)+f(x_0-\delta)-2f(x_0)}{\delta^2}.
\end{equation}
When $f$ is differentiable two times at $x_0$ it is straightforward to see that $f''(x_0)=\overline{D^2f}(x_0)$.
From (\ref{Eq:Bayle}) certainly follows $$I''_M(v)=\overline{D^2I_M}(v)\leq\overline{D^2f_{\bar{v}}}(v)=f_{\bar{v}}''(v),$$ for every $v\in I_{\bar{v}}$.

In a point $\bar{v}$ at which $I_M$ is twice differentiable we observe that $${I}''_M(\bar{v})=\overline{D^2I_M}(\bar{v})\leq f_{\bar{v}}''(\bar{v}).$$ Hence, (\ref{Eq:BP}) yields $$I_M(\bar{v}){I}''_M(\bar{v})\leq I_M(\bar{v})f_{\bar{v}}''(\bar{v})\leq -I_M(\bar{v})\left(\frac{{I'}_M^2(\bar{v})}{n-1}-(n-1)k\right),$$ which is exactly (\ref{GenBP}), because $|II|^2\geq \frac{h^2}{n-1}$, where $h=f_{\bar{v}}'(\bar{v})$ by the first variation formula. If equality holds in (\ref{GenBP}), then $|II|^2= \frac{h^2}{n-1}$, which is equivalent to say that the regular part of $\partial\Omega_{\bar{v}}$ is totally umbilic.  
\end{Dem}
\section{Appendix: Bavard-Pansu}\label{Apendix:BP}
We rewrite for completeness the details of a Theorem that could be immediately deduced from the proof of (i) of \cite{BP} pp. 482, even if that theorem is stated for compact manifolds some of the arguments are still valid for a noncompact manifold satisfying the hypothesis of the theorem below.  
\begin{figure}[h!]
  \centering
    \includegraphics[width=.5\linewidth]{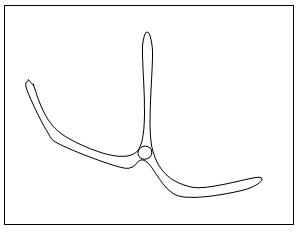}
  \caption{In an isoperimetric region (not in an arbitrary finite perimeter set) of a manifold with Ricci curvature bounded below and non collapsing, it is always possible to put inside a ball of prescribed small volume.}
  \label{fig:fig1.5}
\end{figure} 
\begin{Thm}\label{Thm:BavardPansu}$\cite{BP}$ Let $M^n$ be a complete Riemannian manifold with bounded geometry such that for every volume $v\in ]0, V(M)[$ there exists an isoperimetric region  $\Omega$ of volume $v$. Then $I_M$ is continuous. Moreover  ${I'_M}^+(v), {I_M}^-(v)\leq h=h(n, k,v)$.
\end{Thm}
\begin{Dem}
Let $v\in ]0, V(M)[$ be fixed. Consider a sequence of volumes $v_j\rightarrow v$. By the very definition of the isoperimetric profile we know that $I_M(v_j)\leq A(\partial B_j)$ where $B_j:=B(p,r_j)$ is any geodesic ball inclosing volume $v_j$ and centered at a fixed point $p$. Now take a sequence $\Omega_j$ of isoperimetric regions with $V(\Omega_j)=v_j$, this sequence exists by hypothesis. Theorem 2.1 of \cite{HeintzeKarcher} ensures that the isoperimetric regions have length of mean curvature vector $|H_{\partial\Omega_j}|=:h_j\leq h$, where $h$ is a positive constant that does not depend on $j$ but only on 
$\frac{v}{A(\partial B)}$ where $B$ could be taken as a geodesic ball enclosing volume $v$ in the comparison manifold $\mathbb{M}_k^n$. Again Theorem 2.1 of \cite{HeintzeKarcher} shows that the inradius $\rho_j:=\sup\{d(x, \partial\Omega_j):x\in\Omega_j\}\geq\frac{v}{A(\partial B)}$, if $H_{\partial\Omega_j}$ points inside $\Omega_j$. Observe here that $H_{\partial\Omega_j}$ cannot point outside in the noncompact part if $|H_{\partial\Omega_j}|>1$. If $h_j=|H_{\partial\Omega_j}|\leq 1$ and points outside $\Omega_j$ then $V(\Omega_j)\leq A(\partial\Omega_j)\int_0^{\rho_j}\left(c_k(s)+h_js_k(s)\right)^{n-1}ds$ which implies again that $\rho_j\geq\rho=\rho(n,k, v, A(\partial B))=\rho(n,k,v)>0$. This shows that $\Omega_j$ always contains a geodesic ball of radius $\rho$ centered at a point $p_j$. Now by Theorem \ref{Thm:Uppersemicontinuity} $I_M$ is upper semicontinuous. It remains to show lower semicontinuity. We know that there exists $\bar{v}>0$ such that $V(B_M(q, \rho))\geq\bar{v}>0$ for every $q\in M$, by the noncollapsing hypothesis. Look at the case $v_j\geq v$ then if $v_j-v$ is small enough we can always pick a radius $0<r_j<\rho$ such that $V(B(p_j, r_j))=v_j-v$ again by the noncollapsing hypothesis. Put $\Omega'_j:=\Omega_j\setminus B(p_j, r_j)$, we have $V(\Omega'_j)=v$, thus $I_M(v)\leq A(\partial\Omega')=A(\partial\Omega_j)+A(\partial B(p_j, r_j))$ and finally passing to the limit we obtain $I_M(v)\leq\liminf I_M(v_j)$. If $v_j\leq v$ then the proof is easier and consists in just adding a small ball outside $\Omega_j$ to finish the proof. 
\end{Dem}
\begin{Rem} Applying the proof of Theorem \ref{Thm:BavardPansu} to generalized isoperimetric regions we see easily that the conclusions of Theorem \ref{Thm:BavardPansu} holds if we assume that $M$ has $C^0$-locally bounded geometry.    
\end{Rem}
\begin{Rem} It is not too hard to see that Corollary \ref{BPGen} could be seen also as a corollary of Theorem \ref{Thm:BavardPansu}, without using the proof of Theorem \ref{Main}, because we could argue the continuity of $I_M$ from the proof of Theorem \ref{Thm:BavardPansu} applied to generalized isoperimetric regions and to continue unchanged the proof of Corollary \ref{BPGen}.
\end{Rem}
The argument of the proof of \cite{BP} that cannot be extended easily to the noncompact case with collapsing, concerns the proof of the concavity of the isoperimetric function plus a quadratic function, without passing previously from a proof of the continuity of $I_M$. We don't know if this is possible but a priori the proof seems quite more involved and for the moment we are not able to do it. We present in the following theorem another extension of the arguments of \cite{BP} that permits to argue weaker conclusion on the isoperimetric profile but still not the continuity or concavity.
\begin{Thm}\label{Thm:BavardPansuStrong} Let $M^n$ be a complete Riemannian manifold with $Ricci\geq k$ such that for every volume $v\in ]0, V(M)[$ there exists an isoperimetric region  $\Omega$ of volume $v$. Then for every $[a, b]\subset]0, V(M)[$ there exists a constant $C=C(a, b, n, k, M)$ such that 
$v\mapsto I_M-C(a,b,n,k, M)v^2$ have nonpositive second derivatives in the sense of distributions. 
\end{Thm}
\begin{Dem} If $k<0$ then 
\begin{eqnarray*}
f_{\bar{v}}'' (v) & \leq & -\frac{(n-1)k}{I_M(v)}\\
& \leq & -\frac{(n-1)k}{a}\sup\left\{\frac{\bar{v}}{I_M(\bar{v})}|\bar{v}\in[a, b]\right\}\\
& \leq & -\frac{(n-1)k}{a}\sup\left\{J(h, \rho)|\bar{v}\in[a, b]\right\}\\
& = & -\frac{(n-1)k}{a}\delta\left(n, k, a, b\right),\\
\end{eqnarray*} 
where $J(h,\rho):=\int_0^{\rho}\left((c_k(s)+|h|s_k(s)\right)^{n-1}ds$, $h$ is an upper bound on the length of the mean curvature of the isoperimetric regions in the interval $[a, b]$ and $\rho=\rho(n, k, v, A(\partial B)))=\rho(n,k,v)$, where $B$ is any geodesic ball enclosing a volume $v$ in $\mathbb{M}_k^n$. 
\end{Dem}
\begin{Rem} In our opinion, it remains still an open question whether $Ricci$ bounded below and existence of isoperimetric regions for every volume implies continuity of the isoperimetric profile in presence of collapsing. We are not able to extend to this setting the arguments of $\cite{BP}$, neither to provide a counterexample, because the manifolds with discontinuous isoperimetric profile constructed in $\cite{NardulliPansuDiscontinuous}$ have $Ricci$ curvature tending to $-\infty$.  
\end{Rem}
      \markboth{References}{References}
      \bibliographystyle{alpha}
      \bibliography{ContinuidadeAbraham7}
      \addcontentsline{toc}{section}{\numberline{}References}
      
     \emph{Abraham Mu\~noz Flores\\ Departamento de Matem\'atica\\ Instituto de Matem\'atica\\ UFRJ-Universidade Federal do Rio de Janeiro, Brasil\\ email: abraham@im.ufrj.br\\
     and Departamento de Geometria e Representa\c{c}\~ao Gr\'afica\\
     Instituto de Matem\'atica e Estat\'istica\\
     UERJ-Universidade Estadual do Rio de Janeiro\\
     email: abraham.flores@ime.uerj.br}\\
     
      \emph{Stefano Nardulli\\ Departamento de Matem\'atica\\ Instituto de Matem\'atica\\ UFRJ-Universidade Federal do Rio de Janeiro, Brasil\\ email: nardulli@im.ufrj.br} 
\end{document}